\documentclass[12pt]{article}
\usepackage{amsfonts,amsmath}
\usepackage{amssymb,latexsym}
\usepackage[dvips]{epsfig}
 
\title{Link homology and categorification}
\author{ Mikhail Khovanov\footnotemark[1]}

\date{}

\newtheorem{theorem}{Theorem}

\newcommand{\oplusop}[1]{{\mathop{\oplus}\limits_{#1}}}
 
\begin{document} 

\maketitle
\baselineskip 14pt
 
\def\R{\mathbb R}
\def\Q{\mathbb Q}
\def\Z{\mathbb Z}
\def\C{\mathbb C}
\def\F{\mathbb F} 
\def\o{\otimes}
\def\lra{\longrightarrow}
\def\mc{\mathcal} 
\def\mf{\mathfrak}
\def\lcob{\mathbf{LCob}} 
\def\drawing#1{\begin{center}\epsfig{file=#1}\end{center}}

\begin{abstract} This is a short survey of algebro-combinatorial 
link homology theories which have the Jones polynomial and other 
link polynomials as their Euler characteristics. 
\end{abstract}

{\bf 2000 Mathematics Subject Classification:} 57M25, 57Q45. 

{\bf Keywords and Phrases:} Link homology, quantum link invariants, matrix 
factorizations, Jones polynomial, HOMFLY-PT polynomial. 

\addtocounter{footnote}{1}\footnotetext{ 
Department of Mathematics, Columbia University, 
New York, NY 10027. E-mail: khovanov@math.columbia.edu}

\vspace{0.2in} 

\section{Introduction} 

The discovery of the Jones polynomial by V.~Jones [J] and quantum groups 
by V.~Drinfeld and M.~Jimbo led to an explosive development of 
quantum topology. The newly found topological invariants  
were christened "quantum invariants"; for knots and links they often take 
the form of polynomials. By late 80's-early 90's it was realized that each 
complex simple Lie algebra $\mf{g}$ gives rise to a gaggle of quantum 
invariants. To a link $L$ in $\R^3$ with each component colored by 
an irreducible representation of $\mf{g}$ there is assigned an invariant 
$P(L,\mf{g})$ taking values in the ring of Laurent polynomials 
$\Z[q,q^{-1}]$ (sometimes fractional powers of $q$ are necessary).  
Polynomials $P(L,\mf{g})$ have a representation-theoretical description, via 
intertwiners between tensor products of irreducible representations 
of the quantum group $U_q(\mf{g}),$ the latter a Hopf algebra 
deformation of the universal enveloping algebra of $\mf{g}.$ 
These invariants by no means exhaust all quantum invariants of knots 
and links; various modifications and generalizations include 
finite type (Vassiliev) invariants, invariants associated with quantum 
deformations of Lie superalgebras, etc. 

Quantum $\mf{sl}(n)$ link polynomials, when each component of $L$ 
is colored by the fundamental $n$-dimensional representation, can be conveniently 
encapsulated into a single 2-variable polynomial $P(L),$ known as the 
HOMFLY or HOMFLY-PT polynomial [HOMFLY], [PT]. 

The skein relation 
$$ \lambda P(L_1) - \lambda^{-1} P(L_2) = (q-q^{-1}) P(L_3)$$ 
for any three links $L_1, L_2, L_3$ that differ as shown below

\drawing{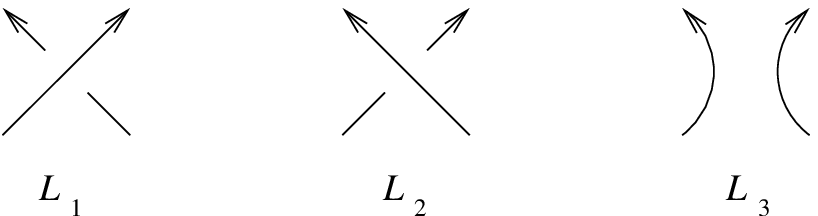}
and the value of $P$ on the unknot, uniquely determines 
the HOMFLY-PT invariant, which lies in the ring $\Z[\lambda^{\pm 1}, (q-q^{-1})^{\pm 1}]$ 
(in the original papers a single variable was used instead of $q-q^{-1},$ making 
$P$ a genuine Laurent polynomial in two variables). 

Specializing $\lambda=q^n,$ for $n \ge 0,$ leads to a link polynomial 
invariant $P_n(L)\in \Z[q,q^{-1}],$ normalized so that  
$P_n(\mathrm{unknot})=q^{n-1}+ q^{n-3} + \dots + q^{1-n}$ for  
$n>0$ and $P_0(\mathrm{unknot})=1.$  

$P_0(L)$ and $P_2(L)$ are the Alexander and Jones polynomials of $L,$ 
respectively, while $P_1(L)$ is a trivial invariant. For $n>0,$ the 
polynomial $P_n(L)$ can be interpreted via the representation theory of 
quantum $\mf{sl}(n),$ and $P_0(L)$ -- via that of the quantum Lie 
superalgebra $U_q(\mf{gl}(1|1)).$ 

The miracle that emerged in the past few years is that these polynomials 
are Euler characteristics of link homology theories: 

\begin{itemize} 
\item The Jones polynomial $P_2(L)$ is the Euler characteristic 
of a bigraded link homology theory $\mc{H}(L),$ discovered in [K1]. 
\item The Alexander polynomial $P_0(L)$ is the Euler characteristic 
of a bigraded knot homology theory, discovered by P.~Ozsv\'ath, Z.~Szab\'o [OS1] and 
J.~Rasmussen [R1]. 
\item The polynomial $P_3(L)$ is the Euler characteristic of 
a link homology theory $\mathrm{H}(L),$ defined in [K2]. 
\item For each $n\ge 1,$ Lev Rozansky and the author constructed 
a bigraded link homology theory $\mathrm{H}_n(L)$ with $P_n(L)$ 
as the Euler characteristic, see [KR1] . 
\item The entire HOMFLY-PT polynomial is the Euler characteristic 
of a triply-graded link homology theory [KR2], [K6] (for a possible alternative 
approach via string theory see [GSV]). 
\end{itemize}

Ideally, a link homology theory should be a monoidal functor $\mc{F}$ 
from the category $\lcob$ of link cobordisms to a tensor triangulated 
category $\mathbf{T}$ (for instance, $\mathbf{T}$ could be the category 
of complexes of $R$-modules, up to chain homotopies, for a 
commutative ring $R$). Objects of $\lcob$ are oriented 
links in $\mathbb{R}^3,$ morphisms from $L_0$ to $L_1$ are isotopy
classes (rel boundary) of oriented surfaces $S$ smoothly and properly embedded 
in $\R^3\times [0,1]$ such that $L_0 \sqcup (-L_1)$ is the boundary of $S$
and $L_i \subset \R^3\times \{i\},$ $i=0,1.$  
In many known examples, $\mc{F}$ is a projective functor: the map  
$\mc{F}(S): \mc{F}(L_0) \lra \mc{F}(L_1)$ 
is well-defined up to overall multiplication by invertible central elements 
of $\mathbf{T}$ (e.g. by $\pm 1$ for homology 
theory $\mc{H}$).

No a priori reason why quantum link invariants should lift 
to link homology theories is known and the general framework 
for lifting quantum invariants to homology theories remains a 
mystery.  We call such a lift a \emph{categorification} of the invariant.
The term categorification was coined by L.~Crane and I.~Frenkel [CF]
in the context of lifting an $n$-dimensional TQFT to an 
$(n+1)$-dimensional one ($n=2,3$ are the main interesting cases).  

Let us also point out that the Casson 
invariant (a degree two finite-type invariant of 3-manifolds) 
is the Euler characteristic of instanton Floer homology, that 
the Seiberg-Witten and Ozsv\'ath-Szab\'o 3-manifold homology 
theories categorify degree one finite-type invariants of 3-manifolds 
(the order of $\mathrm{H}_1(M,\Z)$ when the first homology 
of the 3-manifold $M$ is finite and, more generally, the Alexander polynomial 
of $M$), that equivariant knot signatures are Euler characteristics of 
$\Z/4\Z$-graded link homologies (O.~Collin, B.~Steer [CS], W.~Li), 
and that there exist ideas on how to categorify the 2-variable 
Kauffman polynomial [GW], the colored Jones 
polynomial, and quantum invariants of links colored 
by arbitrary fundamental representations $\Lambda^i V$ 
of $\mf{sl}(n)$ [KR1].

\section {A categorification of the Jones polynomial}
 
In the late nineties the author discovered a homology 
theory $\mc{H}(L)$ of links which is bigraded, 
$$ \mc{H}(L) = \oplusop{i,j\in \Z}\mc{H}^{i,j}(L),$$ 
and has the Jones polynomial as the Euler characteristic, 
$$ P_2(L) = \sum_{i,j\in \Z} (-1)^i q^j \mathrm{rk}(\mc{H}^{i,j}(L)).$$ 
The construction of $\mc{H}$ categorifies the Kauffman bracket 
description of the Jones polynomial. Starting from a  
plane projection $D$ of $L$ we build homology groups $\mc{H}(D)$ 
inductively on the number of crossings of the projection via long exact sequences 

\drawing{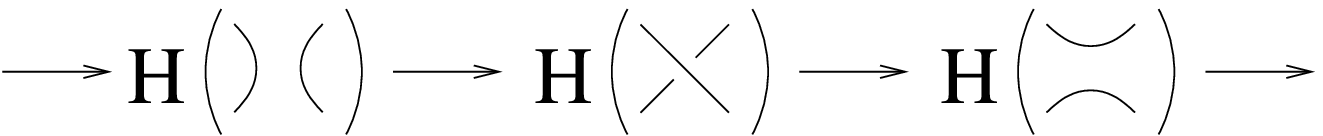} 
and then check that $\mc{H}(D)$ are invariants of $L$ alone. 
Homology of the empty link is $\Z,$ homology of the unknot is 
$\mc{A}=\Z[X]/(X^2),$ which should be thought of 
as the integral cohomology ring of the 2-sphere. 
Homology of the $k$-component unlink is $\mc{A}^{\otimes k}.$ 
The obvious cobordisms between unlinks turn $\mc{A}$ into a 
commutative Frobenius ring, with the trace map  
$\mathrm{tr}(1)=0,$ $\mathrm{tr}(X)=1$ (in any full-fledged 
link homology theory  homology of the unknot is a commutative Frobenius 
algebra over homology of the empty link). 
$\mc{H}(D)$ is the homology of a complex $\mc{C}(D)$ 
constructed in an elementary way from 
direct sums of tensor powers of $\mc{A}$ and the structure maps of 
this Frobenius ring. 

\begin{theorem} There exists a combinatorially defined bigraded homology theory $\mc{H}(L)$ 
of oriented links in $\R^3.$ Groups $\mc{H}^{i,j}(L)$ are finitely-generated and their 
Euler characteristic is the Jones polynomial. 
The theory is functorial: to an oriented cobordism $S$ between 
links $L_0$ and $L_1$ it assigns a homomorphism of groups 
$$ \mc{H}(S): \mc{H}(L_0) \lra \mc{H}(L_1),$$
well-defined up to overall minus sign and 
of bidegree $(0,-\chi(S)),$ where $\chi(S)$ is the Euler characteristic 
of the surface $S.$ 
\end{theorem} 

\noindent 
That $\pm \mc{H}(S)$ is well-defined was proved in [Ja] and [K4] in two different ways. 

\vspace{0.1in} 

Homology theory $\mc{H}$ is manifestly combinatorial and programs 
computing it were written by D.~Bar-Natan, A.~Shumakovitch and J.~Green. 
The earliest program [BN1] led to the conjecture that ranks of
 the homology groups of alternating links are determined 
by the Jones polynomial and the signature. This conjecture was proved by 
E.-S.~Lee [L1]. For arbitrary knots and links, the structure of $\mc{H}$ 
is more complicated than that of the Jones polynomial;  
right now we don't even have a guess at what the rational homology groups
of arbitrary $(n,m)$-torus knots are. 

\vspace{0.06in} 

We next list several interesting applications of $\mc{H}$ and 
related developments.  

\vspace{0.06in} 

\noindent 
1) J.~Rasmussen used $\mc{H}$ and its deformation studied by E.~S.~Lee [L2]
to give a combinatorial proof of the Milnor conjecture that the slice genus 
of the $(p,q)$-torus knot is $\frac{(p-1)(q-1)}{2}$ and of its generalization to all 
positive knots [R2].  This can also be used to show that certain knots are topologically 
but not smoothly slice without having to invoke Donaldson or Seiberg-Witten gauge theories. 
Originally, the Milnor 
conjecture was proved by P.~Kronheimer and T.~Mrowka via the Donaldson theory [KM]. 

\vspace{0.06in} 

\noindent 
2) Lenhard Ng  [N] obtained an upper bound on the Thurston-Bennequin 
number of a Legendrian link from its homology $\mc{H}(L).$ This bound 
is sharp on alternating knots and on all but one or two knots with at most 
10 crossings.
 
\vspace{0.06in} 

\noindent
3) A.~Shumakovitch [S] showed that over the 2-element field homology 
decomposes: $\mc{H}(L, \F_2)\cong 
\widetilde{\mc{H}}(L, \F_2) \otimes \F_2[X]/(X^2),$ where 
$\widetilde{\mc{H}}(L, \F_2)$ is the reduced homology of $L$ with 
coeffcients in $\F_2.$ P.~Ozsv\'ath and Z.~Szab\'o [OS2] discovered a
spectral sequence with the $E^2$-term $\widetilde{\mc{H}}(L, \F_2)$
that converges to the Ozsv\'ath-Szab\'o homology of the double 
branched cover of $L^!.$

\vspace{0.06in} 

\noindent 
4) P.~Seidel and I.~Smith defined a $\Z$-graded homology 
theory of links via Lagrangian intersection Floer homology of 
a certain quiver variety [SS]. Their theory is similar to $\mc{H}$ in many 
respects, and, conjecturally, isomorphic to $\mc{H}$ after the bigrading 
in the latter is collapsed to a single grading. 

\section{Extensions to tangles.} 

The quantum group $U_q(sl(2))$ controls the extension of the Jones 
polynomial to an invariant of tangles, the latter a functor from 
the category of tangles to the category of $U_q(sl(2))$ representations. 
To a tangle $T$ with $n$ bottom and $m$ top endpoints (an $(m,n)$-tangle) 
there is assigned an intertwiner 
$$f(T): V^{\otimes n} \lra V^{\otimes m}$$ 
between tensor powers of the fundamental representation $V$ of $U_q(sl(2)).$ 

A categorification of the invariant $f(T)$ was suggested in [BFK]. We considered 
the category 
$$\mc{O}^n= \oplusop{0\le k \le n} \mc{O}^{k,n-k},$$ 
the direct sum of parabolic subcategories $\mc{O}^{k,n-k}$ 
of a regular block of the highest weight category for $\mf{sl}(n).$ The category 
$\mc{O}^{k,n-k}$ is equivalent to the category of perverse sheaves 
on the Grassmannian of $k$-planes in $\C^n,$ smooth with respect to 
the Schubert stratification. The Grothendieck group of $\mc{O}^n$ is 
naturally isomorphic (after tensoring with $\C$) to $V^{\otimes n},$ 
considered as a representation of $U_{q=1}(sl(2)),$ and  
derived Zuckerman functors in $D^b(\mc{O}^n)$ lift the action of $sl(2)$ 
on $V^{\otimes n}.$ 
We showed that projective functors in $\mc{O}^n$ 
categorify the action of the Temperley-Lieb 
algebra on $V^{\otimes n}$ and conjectured how 
to extend this to arbitrary tangles, by assigning to a tangle $T$ a 
functor $\mc{F}(T)$ between derived categories $D^b(\mc{O}^n)$ 
and $D^b(\mc{O}^m).$ 

Our conjectures were proved by C.~Stroppel [St1], who worked 
with the graded versions $\mc{O}_{gr}^n$ of these categories, 
associated a functor $\mc{F}(T)$ between derived categories 
$D^b(\mc{O}_{gr}^n)$ and $D^b(\mc{O}_{gr}^m)$ to 
each $(m,n)$-tangle $T$ and 
a natural transformation $\mc{F}(S): \mc{F}(T_0)\lra \mc{F}(T_1)$ 
to a tangle cobordism $S$ between tangles $T_0$ and $T_1.$ 
The whole construction is a 2-functor from the 2-category of 
tangle cobordisms to the 2-category whose objects are 
$\C$-linear triangulated categories, 1-morphisms are exact 
functors and 2-morphisms--natural transformations of functors, 
up to rescalings by invertible complex numbers. 
When the tangle is a link $L,$ this theory produces bigraded 
homology groups, conjecturally isomorphic to $\mc{H}(L)\otimes \C.$ 

\vspace{0.06in} 

The braid group action on $V^{\otimes n}$ lifts to a braid group 
action on the derived category $D^b(\mc{O}^n_{gr}).$  
Restricting to the subcategory $D^b(\mc{O}_{gr}^{1,n-1})$ 
results in a categorification of the Burau representation, previously studied in [KS]. 

\vspace{0.06in}

For a more economical extension of the Jones polynomial to tangles, 
we restrict to even tangles (when the number of endpoints on each of the 
two boundary planes is even) and to the subspace of $U_q(sl(2))$-invariants 
$$\mathrm{Inv}(n) = \mathrm{Hom}_{U_q(sl(2))}(\C, V^{\otimes 2n}) $$ 
in $V^{\otimes 2n}.$ The invariant of a $(2m,2n)$-tangle is a linear 
map 
$$f_{inv}(T): \mathrm{Inv}(n) \lra \mathrm{Inv}(m)$$ 
between these subspaces. 

\vspace{0.1in} 

A categorification of $f_{inv}(T)$ was found in [K3], [K4]. We defined a 
graded ring $H^n$ and established an isomorphism 
$$ K(H^n\mathrm{-mod}) \otimes \C \cong \mathrm{Inv}(n)$$ 
between the Grothendieck group (tensored with $\C$) 
of the category of graded finitely-generated $H^n$-modules and 
the space of invariants in $V^{\otimes 2n}.$ To an even  tangle $T$ 
we assigned an exact functor $\mc{T}$ between the derived categories 
of $H^n\mathrm{-mod}$ (this functor induces the map 
$f_{inv}(T)$ on the Grothendieck groups) and to a tangle cobordism--a natural 
transformation of functors. 
This results in a 2-functor from 
the 2-category of cobordisms between even tangles to the 2-category 
of natural transformation between exact functors in triangulated categories. 
Restricting to links, we recover homology groups $\mc{H}(L).$ 
This approach is more elementary than that 
via category $\mc{O},$ and should carry the same amount of 
information. 

The space of invariants $\mathrm{Inv}(n)$ is a subspace of 
$V^{\otimes 2n}(0),$ the weight zero subspace of $V^{\otimes 2n}.$ 
A categorification of this inclusion, recently found by Stroppel [St2],  
 relates rings $H^n$ and parabolic categories $\mc{O}^{n,n}.$ The latter category is 
equivalent to the category of finite-dimensional modules over 
a  $\C$-algebra $A_{n,n},$ explicitly described 
by T.~Braden [B]. There exists a idempotent $e$ in $A_{n,n}$ such 
that $e A_{n,n} e \cong H^n\otimes \C.$ This idempotent 
picks out all self-dual indecomposable projectives in $A_{n,n}.$ 

Rings $H^n$ can also be used to categorify certain level two 
representations of $U_q(sl(m)),$ see [HK]. 

\vspace{0.1in} 

For a more geometric and refined 
approach to invariants of tangles and tangle cobordisms we 
refer the reader to Bar-Natan [BN2]. Some of his generalizations 
of link homology can be thought of as $G$-equivariant 
versions of $\mc{H},$ for various compact subgroups $G$ of 
$SU(2),$ see [K5] for speculations in this direction and for 
an interpretation of the Rasmussen invariant via the 
$SU(2)$-equivariant version of $\mc{H}.$

\section{sl(n) link homology and matrix factorizations} 

\begin{theorem} \label{sln} For each $n>0$ there exists a homology 
theory which associates bigraded homology 
groups
$$ \mathrm{H}_n(L)\cong \oplusop{i,j\in \Z} \mathrm{H}_n^{i,j}(L)$$
to every oriented link in $\R^3.$ The Euler characteristic of $\mathrm{H}_n$ 
is the polynomial invariant $P_n,$ 
$$ P_n(L) = \sum_{i,j\in \Z} (-1)^i q^j \dim_{\Q}(\mathrm{H}_n^{i,j}(L)).$$ 
Homology groups $\mathrm{H}_n^{i,j}(L)$ are finite-dimensional 
$\Q$-vector spaces, and, for a fixed $L,$ only finitely many of them are non-zero. 
This homology is functorial: an oriented link cobordism $S$ between 
$L_0$ and $L_1$ induces a homomorphism 
$$\mathrm{H}_n(S): \mathrm{H}_n(L_0) \lra \mathrm{H}_n(L_1),$$
well-defined up to overall rescaling by nonzero rationals. 
\end{theorem} 

Groups $\mathrm{H}_n(L) $ are constructed in [KR1], where  we start 
with a presentation for $P_n(L)$ as an alternating sum 
\begin{equation} \label{eq1} 
P_n(L) = \sum_{\Gamma} \pm q^{\alpha(\Gamma)} P_n(\Gamma).
\end{equation}  
Here we choose a generic plane projection $D$ of $L$ with $m$ crossings,
 and form the sum over $2^m$ planar trivalent graphs $\Gamma$ which are given  
by replacing each crossing of $D$ by one of the two planar pictures
on the right 
\drawing{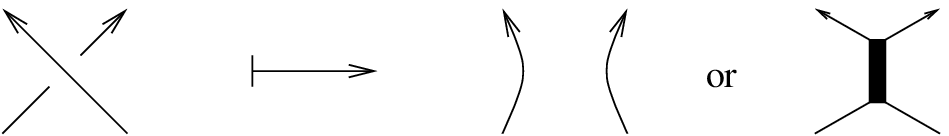} 
Each such planar graph $\Gamma$ has a well-defined invariant 
$P_n(\Gamma) \in \Z[q,q^{-1}],$ with all the coefficients being 
nonnegative integers. Weights $\alpha(\Gamma)$ are given by a simple 
rule. The edges are of two types: regular oriented edges and "wide" unoriented 
edges as on the rightmost picture above.  

We then define single-graded homology groups 
$\mathrm{H}_n(\Gamma)$ which have the graded dimension $P_n(\Gamma)$
and satisfy certain naturality conditions allowing us to build a complex out 
of $\mathrm{H}_n(\Gamma),$ over all modifications $\Gamma$ of the link 
diagram $D.$ The complex is a categorification of the right hand side of the 
equation (\ref{eq1}); its homology groups $\mathrm{H}_n(D)$ 
depend on $L$ only and satisfy the properties listed in Theorem~\ref{sln}. 

Our definition of $\mathrm{H}_n(\Gamma)$ is based on matrix factorizations. 
Let $R=\Q[x_1,\dots, x_k].$ A matrix factorization $M$ of a polynomial 
$f\in R$ consists of a pair of free $R$-modules and 
a pair of $R$-module maps 
 $$ M^0 \stackrel{d}{\lra} M^1 \stackrel{d}{\lra} M^0$$ 
such that $d^2=f\cdot \mathrm{Id}.$ The polynomial $f$ is called 
the \emph{potential} of $M.$ A matrix factorization can be 
thought of as a two-periodic generalized complex; the square of the 
differential in not zero, but a fixed multiple of the identity operator. 
Matrix factorizations were introduced by D.~Eisenbud [E] to study 
homological properties of hypersurface singularities, and later made 
an appearance in string theory, as boundary conditions in Landau-Ginzburg 
models [KL]. 
The tensor product  $M\otimes_R N$ of matrix factorizations with 
potentials $f,g$ is a matrix factorization with potential $f+g.$ 

To each $\Gamma$ we associate a collection of matrix factorizations
$M_1, \dots, M_m,$ one for each crossing of $D,$ with potentials 
$f_1, \dots, f_m$ that add up to zero: $f_1+f_2+\dots + f_m=0.$ 
The tensor product $M_1\otimes M_2 \otimes \dots \otimes M_m$ 
is a two-periodic complex (since the square of the differential is now zero). 
Finally, $\mathrm{H}_n(\Gamma)$ is defined as the cohomology of 
this complex; it inherits a natural $\Z$-grading from that of the 
polynomial algebra $R.$ 

\vspace{0.1in} 

The homology theory $\mathrm{H}_n$ is trivial when $n=1,$ while 
$\mathrm{H}_2(L)\cong \mc{H}(L)\otimes \Q.$ 
The theory $\mathrm{H}_3$ should be closely related to the homology theory constructed 
earlier in [K2] (the two theories have the same Euler characteristic; the one in [K2] 
is defined over $\Z$ and not just over $\Q$). 

J.~Rasmussen [R3] determined homology groups $\mathrm{H}_n(L)$ 
for all 2-bridge knots $L$ and a few other knots (with a mild technical 
restriction $n>4$). Little else is known about homology groups 
$\mathrm{H}_n(L)$ for $n>2.$ 

A Lagrangian intersection Floer homology counterpart of $\mathrm{H}_n$ 
was discovered by C.~Manolescu [M]. His theory $\mc{H}_{(n)symp}(L)$ 
is singly-graded, but defined over $\Z.$ Manolescu conjectured that $\mc{H}_{(n)symp},$ 
after tensoring with $\Q,$ becomes isomorphic to $\mathrm{H}_n,$ with the bigrading 
of the latter folded into a single grading.

\section{Triply-graded link homology and beyond} 

It turns out that the entire HOMFLY-PT polynomial, and not just its 
one-variable specializations, admits a categorification. The 
original construction via degenerate matrix factorizations with a parameter [KR2] 
was later recast in the language of Hochschild homology for bimodules 
over polynomial algebras [K6]. We represent a link $L$ as the closure 
of a braid $\sigma$ with $m$ strands. To $\sigma$ we assign a 
certain complex $F(\sigma)$ of graded bimodules over the polynomial 
algebra $R$ in $m-1$ generators. Taking the Hochschild homology over 
$R$ of each term in the complex produces a complex of bigraded 
vector spaces 
$$ \dots \lra \mathrm{HH}(R, F^j(\sigma)) \lra \mathrm{HH}(R, F^{j+1}(\sigma)) 
  \lra \dots $$
The cohomology groups  $\mathrm{H}(\sigma)$  of this complex 
are triply-graded and depend on $L$ only (a convenient grading normalization 
was pointed out by H.~Wu [W]). The Euler characteristic of 
$\mathrm{H}(L)$ is the HOMFLY-PT polynomial $P(L),$ normalized 
so that $P(\mathrm{unknot})=1.$  

This homology theory suffers from two problems. First, the definition 
requires choosing a braid representative of a knot, rather than just 
a plane projection. Second, it's not possible to assign maps $\mathrm{H}(S)$ 
to all link cobordisms $S$ so as to turn $\mathrm{H}$ into a functor from 
$\lcob$ to the category of (triply-graded) vector spaces (simply because 
homology of the unknot is one-dimensional, while that of an unlink 
is infinite-dimensional). We conjecture that the theory can be redefined 
on $k$-component links for all $k>1$ 
so as to assign finite-dimensional homology groups 
$\widetilde{\mathrm{H}}(L)$ to all oriented links $L,$ and not just 
to knots. The Euler characteristic of $\widetilde{\mathrm{H}}(L)$ will 
still be the HOMFLY-PT polynomial, but rescaled so as to be a Laurent 
polynomial in $\lambda$ and $q$ rather than a rational function. 
The theory should extend to a projective 
functor from the category of \emph{connected} 
link cobordisms to the category of triply-graded vector spaces. 

Further extension of $\widetilde{\mathrm{H}}$
to all link cobordisms should only require a minor modification, where 
one assigns the algebra $\Q[a]$ to the empty link, 
the differential graded algebra 
$$\Q\langle y_1, \dots, y_n\rangle \otimes \Q[a], \hspace{0.1in} 
y_i y_j + y_j y_i =0, \hspace{0.1in} [y_i,a]=0, \hspace{0.1in} 
 d(y_i)=a, \hspace{0.1in} d(a)=0$$ 
to a $k$-component unlink, and suitably resolves each 
$\widetilde{\mathrm{H}}(L),$ viewed as a $\Q[a]$ module with the 
trivial action of $a,$ into a complex of free $\Q[a]$-modules. 

Understanding $\widetilde{\mathrm{H}}$ could be an important step towards 
an algebraic description of knot Floer homology, since we expect 
$\widetilde{\mathrm{H}}$ to degenerate (possibly via a spectral 
sequence) into knot Floer homology of Ozsv\'ath-Szab\'o and 
Rasmussen [OS1], [R1], which categorifies the Alexander polynomial.   

An algebraic description of knot and link Floer homology, if someday found and 
combined with the combinatorial construction [OS3] of Ozsv\'ath-Szab\'o
 3-manifold homology of surgeries on a knot from a filtered version 
of knot Floer homology (and a generalization of their construction to links), 
could lead to a combinatorial definition of Ozsv\'ath-Szab\'o 
and Seiberg-Witten 3-manifold homology and, eventually, 
to an algebraic formulation of gauge-theoretical invariants of 4-manifolds. 
 
\vspace{0.1in} 

In conclusion, we mention two other difficult open problems. 

\vspace{0.1in} 

{\bf I.} Categorify polynomial invariants $P(L,\mf{g})$ of knots 
and links associated to arbitrary complex simple Lie algebras $\mf{g}$ and 
their irreducible representations. 

\vspace{0.1in} 

{\bf II.} Categorify the Witten-Reshetikhin-Turaev invariants of 
3-manifolds. 

\vspace{0.3in}

{\bf References} 

\noindent 
[BN1] D.~Bar-Natan, On Khovanov's categorification of the Jones polynomial, 
\emph{Algebr. Geom. Top.} 2 (2002) 337-370, math.QA/0201043. 

\noindent
[BN2] D.~Bar-Natan, Khovanov's homology for tangles and cobordisms, 
\emph{Geom. Topol.}  9 (2005) 1443-1499,   math.GT/0410495. 

\noindent 
[BFK] J.~Bernstein, I.~B.~Frenkel and M.~Khovanov, A categorification of 
the Temperley-Lieb algebra and Schur quotients of U(sl(2)) via projective 
and Zuckerman functors, Selecta Math. (N.S.) 5 (1999), no. 2, 199--241, 
math.QA/0002087.

\noindent 
[B] T.~Braden, Perverse sheaves on Grassmannians, \emph{Canad. J. Math.} 
54 (2002), no. 3, 493, math.AG/9907152

\noindent 
[CS] O.~Collin and B.~Steer, Instanton Floer homology for knots 
via 3-orbifolds, \emph{J. Differential Geom.} 51 (1), 149-202, 1999.

\noindent 
[CF] L.~Crane, I.~Frenkel, 
Four dimensional topological quantum field theory, 
Hopf categories, and the canonical bases, \emph{Jour. Math. Phys.}
 35 (1994), 5136-5154.

\noindent 
[E] D.~Eisenbud, Homological algebra on a complete intersection, 
with an application to group representations, \emph{Trans. Amer. 
Math. Soc.} {\bf 260} (1980), 35--64. 

\noindent 
[HOMFLY] P.~Freyd, D.~Yetter, J.~Hoste, W.~B.~R.~Lickorish, 
K.~Millett and A.~Ocneanu, A new polynomial invariant of knots 
and links, \emph{Bull. AMS. (N.S.)} {\bf 12} 2, 239--246, 1985. 

\noindent 
[GSV] S.~Gukov, A.~Schwarz and C.~Vafa,
 Khovanov-Rozansky homology and topological strings,
hep-th/0412243. 
 
\noindent 
[GV] S.~Gukov, J.~Walcher, Matrix factorizations and Kauffman 
homology, hep-th/0512298. 

\noindent 
[HK] S.~Huerfano and M.~Khovanov, Categorification of some level 
two representations of sl(n), math.QA/0204333. 

\noindent
[Ja] M.~Jacobsson, An invariant of link cobordisms from Khovanov homology,
\emph{Algebr. Geom. Topol.} 4 (2004) 1211-1251, math.GT/0206303. 

\noindent 
[J] V.~F.~R.~Jones, A polynomial invariant for knots via von Neumann 
algebras, \emph{Bull. Amer. Math. Sco. (N.S.)} {\bf 12} (1985), 103--111. 

\noindent 
[KL] A.~Kapustin and Y.~Li, D-Branes in Landau-Ginzburg models and 
algebraic geometry, hep-th/0210296. 

\noindent 
[K1] M.~Khovanov, A categorification of the Jones polynomial, 
\emph{Duke Math. J.} 101 (2000), no. 3, 359--426, 
math.QA/9908171.

\noindent
[K2] M.~Khovanov, sl(3) link homology,  \emph{Algebr. Geom. Topol.}
 4 (2004) 1045--1081 (electronic), arXiv:math.QA/0304375.

\noindent
[K3] M.~Khovanov, A functor-valued invariant of tangles,
\emph{Algebr. Geom. Topol.} 2 (2002) 665-741,  math.QA/0103190. 

\noindent 
[K4] M.~Khovanov, An invariant of tangle cobordisms,  
\emph{Trans. Amer. Math. Soc.} 358 (2006), 315-327, math.QA/0207264.

\noindent 
[K5] M.~Khovanov, Link homology and Frobenius extensions, 
math.QA/0411447. 

\noindent 
[K6] M.~Khovanov, Triply-graded link homology and Hochschild 
homology of Soergel bimodules, math.GT/0510265. 

\noindent
[KR1] M.~Khovanov and L.~Rozansky, Matrix factorizations and
link homology, math.QA/0401268. 

\noindent
[KR2] M.~Khovanov and L.~Rozansky, Matrix factorizations and
link homology II, math.QA/0505056. 

\noindent 
[KS] M.~Khovanov and P.~Seidel, 
Quivers, Floer cohomology, and braid group actions
\emph{J. Amer. Math. Soc.} 15 (2002), no. 1, 203--271 (electronic), 
math.QA/0006056. 

\noindent 
[KM] P.~Kronheimer and T.~Mrowka,  Gauge theory for 
embedded surfaces I,   \emph{Topology}  32  (1993),  no. 4, 773--826.

\noindent 
[L1] E.~S.~Lee, 
The support of the Khovanov's invariants for alternating knots, 
math.GT/0201105. 

\noindent 
[L2] E.~S.~Lee, An endomorphism of the Khovanov invariant, 
\emph{Adv. Math.} 197 (2), 554--586, (2005). 
math.GT/0210213. 

\noindent 
[N] L.~Ng, A Legendrian Thurston-Bennequin bound from Khovanov homology, 
\emph{Alg. Geom. Topol.} 5 (2005) 1637-1653, 
math.GT/0508649. 

\noindent 
[M] C.~Manolescu, 
Link homology theories from symplectic geometry, math.SG/0601629.

\noindent 
[OS1] P.~Ozsv\'ath, Z.~Szab\'o, Holomorphic disks and knot invariants, 
\emph{Adv. Math.} 186 (1), 58--116 (2004), 
math.GT/0209056. 

\noindent 
[OS2] P.~Ozsv\'ath, Z.~Szab\'o, On the Heegaard Floer homology of 
branched double-covers, math.GT/0309170. 

\noindent 
[OS3] P.~Ozsv\'ath, Z.~Szab\'o, Knot Floer homology and 
rational surgeries, math.GT/0504404. 

\noindent 
[PT] J.~Przytycki and P.~Traczyk, Conway Algebras and
Skein Equivalence of Links, \emph{ Proc. AMS}  {\bf 100}  (1987) 744-748.

\noindent 
[R1] J.~Rasmussen, Floer homology and knot complements, PhD thesis, 
Harvard University, 2003, math.GT/0306378. 

\noindent 
[R2] J.~Rasmussen, Khovanov homology and the slice genus, 
math.GT/0402131. 

\noindent 
[R3] J.~Rasmussen, Khovanov-Rozansky homology of two-bridge 
knots and links, math.GT/0508510. 

\noindent 
[SS] P.~Seidel and I.~Smith, A link invariant from the symplectic 
geometry of nilpotent slices, math.SG/0405089. 

\noindent 
[S] A.~Shumakovitch, Torsion of the Khovanov homology, math.GT/0405474. 

\noindent 
[St1] C.~Stroppel, Categorification of the Temperley-Lieb category, tangles 
and cobordisms via projective functors, \emph{Duke Math. J.} 126 (2005), no. 3,
 547--596.

\noindent 
[St2] C.~Stroppel, Perverse sheaves on Grassmannians, Springer fibres and 
Khovanov homology, math.RT/0608234. 

\noindent 
[W] H.~Wu, Braids, transversal knots and the Khovanov-Rozansky theory, 
math.GT/0508064.

\end{document}